\documentclass[12pt]{article}

\usepackage{amsmath}    
\usepackage{graphicx}   
\usepackage{verbatim}   
\usepackage{color}      
\usepackage{subfigure}  
\usepackage{amssymb}
\usepackage{cite}

\usepackage{tikz}
\usepackage{enumitem}
\def\cp{\,\square\,}

\DeclareMathOperator{\sge}{\rm sg_e}
\DeclareMathOperator{\diam}{\rm diam}

\newtheorem{theorem}{Theorem}[section]
\newtheorem{proposition}[theorem]{Proposition}

\newtheorem{corollary}[theorem]{Corollary}
\newtheorem{lemma}[theorem]{Lemma}

\newcommand{\proof}{\noindent{\bf Proof.\ }}
\newcommand{\qed}{\hfill $\square$\medskip}
\def\cp{\,\square\,}
\usepackage{ifthen}
\newcommand\mreza[2]{
	\pgfmathsetmacro{\n}{#1}
	\pgfmathsetmacro{\m}{#2}
	\ifthenelse{\m=0}{
		\foreach \i in {1,...,\n}{
			\draw (\i,1) -- (\i,2);
			\draw (\i,3.5) -- (\i,4.5);
			\draw[color=white, text=black] (\i,3) node {\tiny \vdots};
		}
		\foreach \j in {1,2,3.5,4.5}{
			\draw (1,\j) -- (\n,\j);
		}
		\foreach \i in {1,...,\n}{
			\draw (\i,3.5) arc [x radius = 0.4em, y radius = 1.7em, start angle = 90, end angle = 270];
			\draw(\i,3.5) arc [x radius = 0.8em, y radius = 2.9em, start angle = 90, end angle = 270];
			\draw(\i,4.5) arc [x radius = 0.8em, y radius = 2.9em, start angle = 90, end angle = 270];
			\draw(\i,4.5) arc [x radius = 1.2em, y radius = 4.1em, start angle = 90, end angle = 270];
			\foreach \j in {1,2,3.5,4.5}{
				\draw[fill=white] (\i,\j) circle(5pt);
		}}
	}
	{
		\ifthenelse{\n = 0}{
			\foreach \i in {1,...,3}{
				\draw (\i,1) -- (\i,\m);
			}
			\draw (5,1) -- (5,\m);
			\draw (6,1) -- (6,\m);
			\foreach \j in {1,...,\m}{
				\draw (1,\j) -- (3.5,\j);
				\draw (4.5,\j) -- (6,\j);
				\draw[color=white, text=black] (4,\j) node {\tiny \dots};
			}
			\foreach \j in {1,...,\m}{
				\draw[fill=white] (1,\j) circle(5pt);
				\draw[fill=white] (2,\j) circle(5pt);
				\draw[fill=white] (3,\j) circle(5pt);
				\draw[fill=white] (5,\j) circle(5pt);
				\draw[fill=white] (6,\j) circle(5pt);
			}
		}{
			\foreach \i in {1,...,\n}{
				\draw (\i,1) -- (\i,\m);
			}
			\foreach \j in {1,...,\m}{
				\draw (1,\j) -- (\n,\j);
			}
			\foreach \i in {1,...,\n}{
				\foreach \j in {1,...,\m}{
					\draw[fill=white] (\i,\j) circle(5pt);
			}}
		}
	}
}

\begin{document}

\title{Strong Edge Geodetic Problem on Complete Multipartite Graphs and some Extremal Graphs for the Problem}
\author{
Sandi Klav\v zar $^{a,b,c}$  \and 	Eva Zmazek $^{a}$
}

\date{\today}

\maketitle

\begin{center}
	$^a$ Faculty of Mathematics and Physics, University of Ljubljana, Slovenia\\
	{\tt sandi.klavzar@fmf.uni-lj.si}\\
	{\tt eva.zmazek@fmf.uni-lj.si}\\
	\medskip
	
	$^b$ Institute of Mathematics, Physics and Mechanics, Ljubljana, Slovenia\\
	\medskip
	
	$^c$ Faculty of Natural Sciences and Mathematics, University of Maribor, Slovenia\\
\end{center}

\begin{abstract}
A set of vertices $X$ of a graph $G$ is a strong edge geodetic set if to any pair of vertices from $X$ we can assign one (or zero) shortest path between them such that every edge of $G$ is contained in at least one on these paths. The cardinality of a smallest strong edge geodetic set of $G$ is the strong edge geodetic number ${\rm sg_e}(G)$ of $G$. In this paper, the strong edge geodetic number of complete multipartite graphs is determined. Graphs $G$ with ${\rm sg_e}(G) = n(G)$ are characterized and ${\rm sg_e}$ is determined for Cartesian products $P_n\,\square\, K_m$. The latter result in particular corrects an error from the literature. 
\end{abstract}


\medskip
\noindent{\bf AMS Subj.\ Class.: 05C12, 05C70}

\section{Introduction}

Covering vertices or edges of a graph by the smallest number of paths is a fundamental optimization problem and appears in the literature in several variations depending on the properties one requires from the paths. In the isometric path cover problem (alias geodetic cover problem) the aim is to cover all the vertices by a minimum number of shortest paths~\cite{chakra-2023a+, chakra-2023b+, Fitz1999, Fitz2001, Manuel19, pan-2006}. In the path cover problem we want to cover all the vertices by a minimum number of vertex disjoint paths~\cite{chen-2022, custic-2022, montarezi-2020}. Dual concepts have also been studied as for instance the $k$-path covers which are sets $S$ of vertices of a graph $G$ such that every path of order $k$ in $G$ contains at least one vertex from $S$, see~\cite{brause-2017, bresar-2013, erves-2022}. In the edge version of the isometric path cover problem we want to cover all the edges by a minimum number of shortest paths~\cite{anand-2018, san-2017, san-2020}. In this paper we are interested in the strong edge geodetic problem introduced in~\cite{manuel-2017} as follows. 

Let $G = (V(G), E(G))$ be a graph. A set of vertices $X\subseteq V(G)$ is a {\em strong edge geodetic set} if to any pair of vertices $u$ and $v$ from $X$ we can assign a shortest $u,v$-path $P_{uv}$ such that every edge $xy\in E(G)$ is contained in at least one on the paths $P_{uv}$. The cardinality of a smallest strong edge geodetic set of $G$ is the {\em strong edge geodetic number} $\sge(G)$ of $G$. Such a set is briefly called  a {\em $\sge(G)$-set}. 

In the seminal paper~\cite{manuel-2017}  it was proved, among other results, that the strong edge geodetic problem is ${\cal NP}$-complete. In~\cite{davot-2021} it was further proved that there is no approximation of the strong edge geodetic number with an approximation factor better that $781/780$. 
Several additional results on the strong edge geodetic number were reported in~\cite{xavier-2020, zmazek-2021}. In the latter paper, the strong edge geodetic number was determined for Cartesian products $P_n\cp P_k$, where $k\in \{2,3,4\}$.  

The vertex version of the strong edge geodetic problem is known as the {\em strong geodetic problem} and was studied for the first time in~\cite{manuel-2020, manuel-2018}. The strong geodetic problem is also ${\cal NP}$-complete and remains such even when restricted to bipartite graphs and multipartite graphs~\cite{irsic-2018}. Moreover, determining whether a given set $X$ is a strong geodetic set is ${\cal NP}$-hard~\cite{davot-2021} as well. 

The strong geodetic number of complete bipartite (resp.\ mutipartite) graphs received a lot of attention. First, in~\cite{irsic-2018} the problem was solved for balanced complete bipartite graphs $K_{n,n}$. Subsequently, using different approaches, a formula for the strong geodetic number of arbitrary complete bipartite graphs was derived in~\cite{irsic-2019} and in~\cite{gledel-2020}. In~\cite{irsic-2019}, a lower bound for the strong geodetic number of a complete multipartite graph was given and it was conjectured that the strong geodetic number remains ${\cal NP}$-complete on complete mutipartite graphs. In~\cite{davot-2021} this conjecture was disproved by developing a polynomial algorithm for the strong geodetic number of complete mutipartite graphs. In this direction we emphasize that in~\cite{mezzini-2022} an $O(|E(G)|\cdot |V(G)|^2)$ algorithm for computing the strong geodetic number of an outerplanar graph $G$ was developed. Several additional interesting results on the strong geodetic problem were presented in~\cite{wang-2020}. Among other results, relations between the strong geodetic number and the connectivity and the diameter were established, and graphs with the strong geodetic number equal to $2$, $|V(G)|-1$, and $|V(G)|$ were characterized. 

Motivated by the efforts to determine the strong geodetic number of complete bipartite graphs, we determine in Section~\ref{sec:multipartite} the strong edge geodetic number of complete bipartite graphs and using this result we then determine the strong edge geodetic number of arbitrary complete multipartite graphs. In Section~\ref{sec:largest} we characterize graphs $G$ with $\sge(G) = n(G)$ and discuss the graphs with $\sge(G) = n(G)-1$. In particular we observe that Cartesian products $P_2\cp K_n$ belong to this family of graphs. This corrects~\cite[Theorem 13]{xavier-2020} where it is wrongly stated that $\sge(P_2\cp K_n) = 2n-2$. We then proceed by determining $\sge(P_m\cp K_n)$ for all $m, n\ge 2$. 

We conclude the introduction by giving some definitions needed. The order of a graph $G$ is denoted by $n(G)$. A vertex $u$ of a graph $G$ is {\em universal} if $\deg_G(u) = n(G) - 1$. The {\em Cartesian product} $G\cp H$ of graphs $G$ and $H$ is the graph with the vertex set $V(G) \times V(H)$, vertices $(g,h)$ and $(g',h')$ being adjacent if either $g=g'$ and $hh'\in E(H)$, or $h=h'$ and $gg'\in E(G)$. As usual, $\chi'(G)$ is the chromatic index of $G$. For a positive integer $n$, the set $\{1,\ldots, n\}$ will be dented by $[n]$.  

If $U$ is a strong edge geodetic set, then we will denote by $\widehat{U}$ the set of associated paths with endpoints from $U$ which cover all the edges of $G$. Clearly, $\widehat{U}$ is not unique, but unless stated otherwise, we will assume that $\widehat{U}$ has been selected and is fixed.

\section{Complete multipartite graphs}
\label{sec:multipartite}

In this section we determine the strong edge geodetic number of complete multipartite graphs. To do so, we first prove the corresponding result for complete bipartite graphs which reads as follows. 


\begin{theorem}
\label{thm:complete bipartite}
	If $n \geq m \geq 2$, then the following hold.
	\begin{enumerate}[label=(\roman*)]
		\item If $n$ is even, then
		$$
		\sge(K_{n,m}) =
		\begin{cases}
			n+1; & n=m, \\
			n; & n \geq m+1.
		\end{cases}
		$$
		\item If $n$ is odd, then
		$$
		\sge(K_{n,m}) =
		\begin{cases}
			n+2; & n=m, \\
			n+1; & n=m+1, \\
			n; & n \geq m+2.
		\end{cases}
		$$
	\end{enumerate}
\end{theorem}

In the rest of the section we assume that $n \geq m \geq 2$ and that the bipartition of $K_{n,m}$ is $(X,Y)$, where $X=\{x_0, \ldots, x_{n-1}\}$ and $Y=\{y_0, \ldots, y_{m-1}\}$. 

\begin{lemma}
\label{lemma: atleastonepartofKnm}
	If $U$ is a strong edge geodetic set of $K_{n,m}$, then $X \subseteq U$ or $Y \subseteq U$.
\end{lemma}

\proof
	Let $U$ be a strong edge geodetic set of the graph $K_{n,m}$. Suppose on the contrary that there exist vertices $x_i \in X \setminus U$ and $y_j \in Y \setminus U$. Because $\diam(K_{n,m})=2$ and $x_iy_j$ is an edge of $K_{n,m}$, none of the shortest paths with endpoints from $U$ can cover the edge $x_iy_j$, that is, $U$ cannot be a strong edge geodetic set.
\qed

\begin{lemma}
\label{lem: atleastn+1}
If $U$ is a strong edge geodetic set of $K_{n,m}$ and $Y \subseteq U$, then $|U| \geq n+1$.
\end{lemma}

\proof
	Suppose $U$ is a strong edge geodetic set of $K_{n,m}$, where $U = Y \cup X'$ with $X' \subseteq X$ and $|X'| = k$, $0 \leq k \leq n$. Consider an arbitrary vertex $y_j \in Y$. There are exactly $n-k$ edges between $y_j$ and $X \setminus X'$. Because the shortest paths that cover these edges have both of their endpoints in $Y$, it has to hold $m-1 \geq n-k$. This in turn implies that $|U| = |Y| + |X'| = m+k \geq n+1$.
\qed

\begin{corollary}
\label{cor:lower-even}
If $n \geq m \geq 2$, then $\sge(K_{n,m}) \geq n$. Moreover, if $m=n$, then $\sge(K_{n,n}) \geq n+1$.
\end{corollary}

\proof
If $m=n$, then the second assertion of the corollary follows immediately from Lemmas~\ref{lemma: atleastonepartofKnm} and~\ref{lem: atleastn+1}. Suppose now that $n> m$ and let $U$ be a smallest strong edge geodetic set of $K_{n,m}$, so that $|U| = \sge(K_{n,m})$. Then $X \subseteq U$ or $Y \subseteq U$ by Lemma~\ref{lemma: atleastonepartofKnm}. If $X \subseteq U$, then $\sge(K_{n,m}) = |U| \ge |X| = n$. And if $Y \subseteq U$, then $\sge(K_{n,m}) \ge n +1$ follows by  Lemma~\ref{lem: atleastn+1}. 
\qed

We have thus established the lower bound for the case when $n$ is even. For $n$ odd we proceed as follows. 

\begin{lemma}
\label{lem:odd}
Let $U$ be a strong edge geodetic set of $K_{n,m}$. If $n$ is odd and $X \subseteq U$, then $|U| \geq \frac{2n}{n+1} + m$.
\end{lemma}

\proof
Let $U$ be  a strong edge geodetic set of $K_{n,m}$, where $U=X \cup Y'$ with $Y' \subseteq Y$ and $|Y'| = k$. For each edge $xy$, where $x \in X$ and $y \in Y'$, we put the shortest path $xy$ to $\widehat{U}$.
The edges between vertices from $X$ and $Y \setminus Y'$ must be covered by the shortest paths of length $2$ with both of their endpoints in $X$. For each pair of vertices from $X$ we can put only one shortest path to $\widehat{U}$, so we can only put $n \choose 2$ shortest paths to $\widehat{U}$ to cover the $n \cdot (m-k)$ edges between the vertices from $X$ and the vertices from $Y \setminus Y'$. Moreover, because the degree of every vertex from $Y \setminus Y'$ is $n$, which we have assumed to be odd, each vertex from $Y \setminus Y'$ must be the central vertex of at least $(n+1)/2$ shortest paths from $\widehat{U}$. Since $U$ is a strong edge geodetic set this implies that ${n \choose 2} \geq (m-k) \cdot \frac{n+1}{2}$. This inequality rewrites to $k \ge m - n(n-1)/(n+1)$ which in turn implies that $|U| = n+k \geq n+m- n(n-1)/(n+1) = \frac{2n}{n+1} + m$.
\qed

\begin{corollary}
	If $n\ge 3$ is odd, then $\sge(K_{n,n}) \geq n+2$ and $\sge(K_{n,n-1}) \geq n+1$. 
\end{corollary}

\proof
Let $U$ be a smallest strong edge geodetic set of $K_{n,n}$. By Lemma~\ref{lem:odd}, 
$$|U| \geq  n + \frac{2n}{n+1}\,.$$
As $|U|$ is an integer and $2n/(n+1) > 1$ for $n\ge 2$ we get  $|U| = \sge(K_{n,n}) \ge n + 2$. 

Let now $U$ be a smallest strong edge geodetic set of $K_{n,n-1}$. By Lemma~\ref{lemma: atleastonepartofKnm} we have $X \subseteq U$ or $Y \subseteq U$. In the latter case,   Lemma~\ref{lem: atleastn+1} gives $\sge(K_{n,n-1}) \geq n+1$. Assume second that $X \subseteq U$. Then Lemma~\ref{lem:odd} gives
$$|U| \geq \frac{2n}{n+1} + (n-1) = n + \frac{n-1}{n+1}\,.$$
Since  $\frac{n-1}{n+1} > 0$ for $n \geq 2$ and since $|U|$ is an integer, also in this case we get $\sge(K_{n,n-1}) \geq n+1$.
\qed

So far, we have proved the lower bound for all the cases of Theorem~\ref{thm:complete bipartite}. In the following we will construct in each case a strong edge geodetic set of the required cardinality. 

\medskip\noindent
{\bf Case 1}: $n$ is even. \\
We first consider $K_{n,n}$ and prove that 
\begin{equation}
\label{eq:n-n-even}
\sge(K_{n,n}) \leq n+1\,.
\end{equation}
We claim that $U = X \cup \{y_{n-1}\}$ is a strong edge geodetic set of $K_{n,n}$. For every $0\le i\le n-1$, add the shortest path $x_i y_{n-1}$ to $\widehat{U}$ to cover the edge $x_i y_{n-1}$. Then all the other edges must be covered by shortest paths of the form $x_i y_j x_k$, where $i \not=k$. To do so, we use edge-colorings of $K_n$. It is well-known that $\chi'(K_n) = n-1$ for even $n$. Let $V(K_n) = \{0,1,\dots,n-1\}$. Then an edge-coloring $c$ of $K_n$ using $n-1$ colors can be defined as follows: if $i,j \in \{0,1,\dots,n-2\}$, $i \not= j$, then let $c(ij) = (i+j) \mod (n-1)$, and for for $i \in \{0,1,\dots,n-2\}$ let $c(i(n-1)) = 2i \mod (n-1)$.

\begin{figure}[ht!]
\begin{center}
\begin{tikzpicture}[scale=1.0,style=thick]
\def\vr{3pt}

\begin{scope}[xshift=0cm, yshift=0cm]
\path (0,0) coordinate (x1);
\path (2,0) coordinate (x2);
\path (4,0) coordinate (x3);
\path (6,0) coordinate (x4);
\path (8,0) coordinate (x5);
\path (10,0) coordinate (x6);
\path (0,3) coordinate (y1);
\path (2,3) coordinate (y2);
\path (4,3) coordinate (y3);
\path (6,3) coordinate (y4);
\path (8,3) coordinate (y5);
\path (10,3) coordinate (y6);
\foreach \i in {1,...,6}
{ \draw (y1) -- (x\i); }
\foreach \i in {1,...,6}
{ \draw (y2) -- (x\i); }
\foreach \i in {1,...,6}
{ \draw (y4) -- (x\i); }
\foreach \i in {1,...,6}
{ \draw (y5) -- (x\i); }
\foreach \i in {1,...,6}
{ \draw (y6) -- (x\i); }
\draw[densely dotted] (x1) -- (y3) -- (x3);
\draw[loosely dashed] (x4) -- (y3) -- (x5);
\draw[line width=0.8mm] (x2) -- (y3) -- (x6); 
\foreach \i in {1,...,6}
{
\draw (x\i)  [fill=white] circle (\vr);
\draw (y\i)  [fill=white] circle (\vr);
}
\foreach \i in {1,...,6}
{
\draw (x\i)  [fill=black] circle (\vr);
}
\draw (y6)  [fill=black] circle (\vr);
\draw[below] (x1) node {$x_0$}; 
\draw[below] (x2) node {$x_1$}; 
\draw[below] (x3) node {$x_2$}; 
\draw[below] (x4) node {$x_3$}; 
\draw[below] (x5) node {$x_4$}; 
\draw[below] (x6) node {$x_5$}; 
\draw[above] (y3) node {$y_2$}; 
\end{scope}

\end{tikzpicture}
\end{center}
\caption{Shortest paths from $\widehat{U}$ that cover edges incident to $y_2$.}
\label{fig:four-paths}
\end{figure}
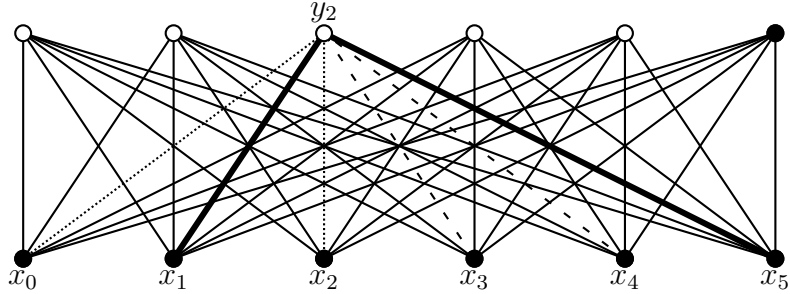

In the covering of $K_{n,n}$ that we are constructing, we put the shortest path $x_i y_j x_k$ to $\widehat{U}$ if and only if $c(ik) = j$. See Fig.~\ref{fig:four-paths}, where this construction is illustrated for the case $n=6$ and the edges incident to $y_2$. In $K_6$, we have $c(02) = c(15) = c(34) = 2$, hence the paths $x_0y_2x_2$, $x_1y_2x_5$, and $x_3y_2x_4$ belong to $\widehat{U}$. 

Using this construction, a pair of vertices $x_i$ and $x_k$ is never used twice, and for each vertex $y \in Y \setminus \{y_{n-1}\}$ the shortest paths in $\widehat{U}$ have pairwise different endpoints. Since in $c$ every color is used exactly $n/2$ times, the shortest paths from $\widehat{U}$ passing through $y_j$ cover all the edges incident with $y_j$. This proves~\eqref{eq:n-n-even}. 

Consider now $K_{n,m}$, where $m \leq n-1$ (and $n$ is even). We need to show that  $\sge(K_{n,m}) \leq n$. For this sake we claim that $X$ is a strong edge geodetic set. Indeed, use the above edge-coloring $c$ of $K_n$ and for each $y_i \in Y$, $i\in \{0, \ldots, m-1\}$, put  all the shortest paths $x_j y_i x_k$ to $\widehat{U}$ for which $c(jk) = i$. By the above argument, $X$ is indeed a strong edge geodetic set and hence $\sge(K_{n,m}) \leq n$ in this subcase. 

\medskip\noindent
{\bf Case 2}: $n$ is odd. \\
We first consider $K_{n,n}$ and prove that $\sge(K_{n,n}) \leq n+2$. For this purpose consider the set $U = X \cup \{y_{n-2}, y_{n-1}\}$. The subgraph of $K_{n,n}$ induced by the set of vertices $V(K_{n,n}) \setminus \{x_{n-1}, y_{n-1}\}$ is isomorphic to $K_{n-1,n-1}$. As $n-1$ is even, we can cover its edges by the paths as described in Case 1 to derive~\eqref{eq:n-n-even}. Recall that for this covering, the vertices $x_0, \ldots, x_{n-1}$ and $y_{n-2}$ are used. To cover the edges $y_{n-1}x_i$, $0\le i\le n-1$, add the shortest paths $y_{n-1}x_i$ to $\widehat{U}$. Finally, to cover the remaining yet uncovered edges, that is, the edges $x_{n-1} y_i$, where $i \in \{0, \dots, n-2\}$, put the shortest paths $x_{n-1} y_i x_i$ to $\widehat{U}$.

We next show that $\sge(K_{n,n-1}) \leq n+1$. In this subcase set $U = X \cup \{y_{n-2}\}$. Then as in the above subcase, cover the edges of the subgraph of $K_{n,n-1}$ induced by the set of vertices $V(K_{n,n-1}) \setminus \{x_{n-1}\}$ as described in Case 1 to derive~\eqref{eq:n-n-even}. After that, to cover the edges $x_{n-1} y_i$, where $i \in \{0, \dots, n-2\}$ we add to $\widehat{U}$ the shortest paths $x_{n-1} y_i x_i$.

Consider finally $K_{n,m}$, where $m \leq n-2$. In this case, $X$ is a  strong edge geodetic set. For this sake note that by the second subcase of Case 1 we know that $\{x_0, \ldots, x_{n-2}\}$ is a strong edge geodetic set of the subgraph of $K_{n,m}$ induced by the set $V(K_{n,m})\setminus \{x_{n-1}\}$. To cover the remaining not yet covered edges $x_{n-1} y_i$, where $i \in \{0, \dots, m-1\}$ we add to $\widehat{U}$ the shortest paths $x_{n-1} y_i x_i$. From here it is clear that $X$ is a  strong edge geodetic set of $K_{n,m}$ and we conclude that in this subcase $\sge(K_{n,m}) \le n$. 

We have thus established all the upper bounds which completes the proof of Theorem~\ref{thm:complete bipartite}. Using it we can in turn determine the strong edge geodetic number of complete multipartite graphs as follows. 

\begin{theorem}
\label{thm:complete multipartite}
	If $k\ge 2$ and $2\le n_1\le n_2\ \le \cdots \le n_k$, then the following hold.
	\begin{enumerate}[label=(\roman*)]
		\item If $n_1$ is even, then
		$$
		\sge(K_{n_1,\dots, n_k}) =
		\begin{cases}
			\sum\limits_{j=2}^{k} n_j + 1; & n_2\in \{n_1, n_1+1\}, \\
			\sum\limits_{j=2}^{k} n_j; & \mbox{otherwise}, 
		\end{cases}
		$$
		\item If $n_1$ is odd, then
		$$
		\sge(K_{n_1,\dots, n_k}) =
		\begin{cases}
			\sum\limits_{j=2}^{k} n_j + 2; & n_2 = n_1, \\
			\sum\limits_{j=2}^{k} n_j; & \mbox{otherwise}, 
		\end{cases}
		$$
	\end{enumerate}
\end{theorem}

\proof
Let $k\ge 2$ and $2\le n_1\le n_2\ \le \cdots \le n_k$, and set $G = K_{n_1,\dots, n_k}$ for the rest of the proof. Let $X_i$, $i\in [k]$, be the partition sets of $G$, where $|X_i| = n_i$. Let $U$ be an arbitrary (smallest) strong edge geodetic set of $G$. If $i\ne j$, then, by Lemma~\ref{lemma: atleastonepartofKnm}, we see that $X_i \subseteq U$ or $X_j\subseteq U$. If follows that $U$ contains $k-1$ of the partite sets. 

Let $W = \cup_{i=2}^{k}X_i$. Then it remains to cover the edges between $X_1$ and each of the $X_i$, $i\ge 2$. More precisely, we need to cover the edges in induced subgraphs $K_{n_1,n_i}$, $i\ge 2$, where the partite sets of cardinality $n_i$ are already included. Clearly, no shortest path of $G$ has length greater than $2$ and a shortest path of length $2$ has both endpoints in the same set of the partition. Hence an edge of $K_{n_1, n_i}$ can be covered only by vertices in $X_1 \cup X_i$, for every $2 \le i \le k$. 

Assume first that $n_2 \ge n_1 + 2$. Then by Theorem~\ref{thm:complete bipartite} and its proof we infer that $X_2$ is a strong edge geodetic set of $K_{n_1,n_2}$ with minimum cardinality. Furthermore, since if $k > 2$, for every $2 < i \le k$ we have $n_i \ge n_2 \ge n_1 + 2$, again by Theorem~\ref{thm:complete bipartite} we have that $X_i$ is a strong edge geodetic set of $K_{n_1, n_i}$ with minimum cardinality. Therefore, $W$ is a is a strong edge geodetic set of $G$ with minimum cardinality no matter whether $n_1$ is even or odd. Moreover, we get the same conclusion if $n_1$ is odd and $n_2 = n_1 + 1$. Assume next that $n_1$ is odd and $n_2 = n_1$. Then Theorem~\ref{thm:complete bipartite}(ii) implies that the union of $X_2$ and two vertices of $X_1$, say $u$ and $w$ is a strong edge geodetic set of $K_{n_1,n_2}$ with minimum cardinality. In this case we conclude that $W\cup \{u,v\}$ is a strong edge geodetic set of $G$ with minimum cardinality. The cases when $n_1$ is even and $n_2\in \{n_1, n_1+1\}$ are treated similarly. 
\qed

\section{Graph with large strong edge geodetic sets}
\label{sec:largest}

In this section we first characterize graphs $G$ with $\sge(G) = n(G)$. After that we consider graphs $G$ with $\sge(G) = n(G) - 1$ and determine $\sge(P_n\cp K_m)$. In particular, $\sge(P_2\cp K_m) = 2m - 1$, which corrects a result from~\cite{xavier-2020}. 

Let $G$ be a graph and $uv\in E(G)$. We say that a vertex $v$ is a {\em dominant neighbor} of $u$ if $N[u]\subseteq N[v]$, where $N[u] = \{u\} \cup \{x:\ ux\in E(G)\}$ is the {\em closed neighborhood} of a vertex $u$. Vertices $u$ and $v$ of a graph $G$ are {\em twins} if $N[u] = N[v]$. Note that twins are necessarily adjacent and that if $u$ and $v$ are twins, then $u$ is a dominant neighbor of $v$ and $v$ is a dominant neighbor of $u$.  

The following lemma seems to be of independent interest. 

\begin{lemma}
\label{lem:useful}
Let $G$ be a graph and $U\subseteq V(G)$ a strong edge geodetic set. If $v$ is a dominant neighbor of $u$, then $u\in U$. In particular, if $u$ and $v$ are twin vertices, then $u\in U$ and $v\in U$. 
\end{lemma}

\proof
Let $uv\in E(G)$ and $N[u]\subseteq N[v]$. If $P$ is a shortest path in $G$ which contains the edge $uv$, then one of the end-points of $P$ must be $u$, for otherwise $P$ would not be shortest. If further $u$ and $v$ are twins, then also $N[v]\subseteq N[u]$ and thus also $v\in U$.  
\qed

\begin{proposition}
\label{prop:sge=n}
Let $G$ be a graph. Then $\sge(G) = n(G)$ if and only if every vertex of $G$ has a dominant neighbor.  
\end{proposition}

\proof
If every vertex of $G$ has a dominant neighbor, then every vertex lies in every strong edge geodetic set by Lemma~\ref{lem:useful}. Hence $\sge(G) = n(G)$. 

Assume now that a vertex $u\in V(G)$ does not admit a dominant neighbor. We claim that $U = V(G)\setminus \{u\}$ is a strong edge geodetic set of $G$. Let $v$ be an arbitrary neighbor of $u$. Since $N[u] \not \subseteq N[v]$, there exists a vertex $w\in N[u]\setminus N[v]$. To cover the edge $uv$, put the shortest path $wuv$ to $\widehat{U}$. Proceed analogously for every neigbor $v'$ of $u$, where if the edge $v'u$ has been already covered before, do nothing. In this way all edges incident with $u$ are covered. Let next $xy$ be an arbitrary edge from $E(G)$ where $\{x, y\} \cap \{u\} = \emptyset$. Then add to $\widehat{U}$ the shortest path $xy$. Clearly, the paths added so far to $\widehat{U}$ cover all the edges of $G$ and we conclude that $\sge(G) <  n(G)$. 
\qed

Proposition~\ref{prop:sge=n} implies several results from~\cite{xavier-2020} as for instance~\cite[Theorem~8]{xavier-2020} which asserts that if a graph $G$ contains at least two universal vertices, then $\sge(G) = n(G)$. 

A vertex $u$ of a graph $G$ is {\em simplicial} if $N(u)$ induces a clique of $G$. If $u$ is a simplicial vertex and $v$ its arbitrary neighbor, then $N[u] \subseteq N[v]$. Denoting by $s(G)$ the number of simplicial vertices of $G$ Lemma~\ref{lem:useful} thus implies:  

\begin{corollary}
\label{cor:simplicial}
If $G$ is a graph, then $\sge(G) \ge s(G)$. 
\end{corollary}

Since in $K_n$ every vertex is simplicial, Corollary~\ref{cor:simplicial} implies that $\sge(K_n) = n$. We can also deduce this fact from Lemma~\ref{lem:useful} by observing that each pair of vertices of $K_n$ are twins. 

Lemma~\ref{lem:useful} implies also the following. 

\begin{corollary}
\label{cor:n-1}
	If a graph $G$ contains a universal vertex, then $\sge(G) \ge n(G)-1$. Moreover, if there is only one universal vertex, then $\sge(G) = n(G)-1$.
\end{corollary}

\proof
Let $w$ be a universal vertex of $G$. Then $w$ is a dominant neighbor of every vertex $u \in V(G) \setminus \{w\}$, hence Lemma~\ref{lem:useful} implies that $V(G) \setminus \{w\} \subseteq U$ for every strong edge geodetic set $U$ of $G$. Thus $\sge(G) \geq n(G)-1$. In the case when $w$ is a unique universal vertex of $G$, then with the same arguments as we had in the last part of the proof of Proposition~\ref{prop:sge=n} we infer that $V(G)\setminus \{w\}$ is a strong edge geodetic set. Hence $\sge(G) \leq n(G)-1$ when $G$ has a unique universal vertex, so that in this case $\sge(G) = n(G)-1$. 
\qed

The second assertion of Corollary~\ref{cor:n-1} was earlier presented as~\cite[Theorem~5]{xavier-2020}. Moreover, in~\cite[Theorem~5]{xavier-2020} it was also claimed that $\sge(P_2 \cp K_m) = 2m - 2$. It can be checked that the result is not true and that instead the Cartesian products $P_2 \cp K_m$ also belong to the family of graphs $G$ for which $\sge(G) \ge n(G)-1$. More generally, we have the following result. 

\begin{theorem}
\label{thm:path-by-complete}
	If $m \geq 3$ and $n \geq 2$, then
	$$\sge(P_n \cp K_m) =
	\begin{cases}
		m k; & n =k^2, \\
		m k + (m-1); & n = k^2 + h, 1 \leq h \leq k, \\
		m k + m; & n = k^2+h, k+1 \leq h \leq 2k.
	\end{cases}
	$$
\end{theorem}

\proof
Set $V(K_m) = [m]$ and $V(P_n) = [n]$ where $i(i+1) \in E(P_n)$ for $i \in [n-1]$. If $y\in V(K_m)$, then we will denote by $P_n^y$ the subgraph of $P_n\cp K_m$ induced by the vertices $(i,y)$, $i\in [n]$. $P_n^y$ is also called a $P_n$-layer of $P_n\cp K_m$ and is isomorphic to $P_n$. Throughout the proof we will use the fact that in a shortest path of $P_n \cp K_m$ there is at most one edge between two distinct $P_n$-layers. 

Consider first the case $n=k^2$, where $k \in \mathbb{N}$. In this case we claim that $U_1= \bigcup_{i=1}^k \bigcup_{j=1}^m \{(i^2,j)\}$ is a strong edge geodetic set of $P_n \cp K_m$. To cover all the edges of $P_n \cp K_m$ we proceed as follows. For every $j \in [m]$, put to $\widehat{U_1}$ the unique shortest path between the vertices $(1,j)$ and $(k^2,j)$. For every pair $y_1,y_2 \in V(K_m)$, $y_1 < y_2$, we put the following shortest paths to $\widehat{U_1}$. 
\begin{itemize}
\item For every $i\in [k]$, put to $\widehat{U_1}$ the unique shortest path (of length $1$) between the vertices $(i^2,y_1)$ and $(i^2,y_2)$. 
\item For every $2\le i \le k$, and for every $l \in [i-1]$, put to $\widehat{U_1}$ the shortest path between the vertices $(l^2,y_1)$ and $(i^2,y_2)$ that contains the edge $((i-1)^2+l,y_1)((i-1)^2+l,y_2)$, and the shortest path between the vertices $(l^2,y_2)$ and $(i^2,y_1)$ that passes through the edge $(i(i-1)+l,y_1)(i(i-1)+l,y_2)$. See Fig.~\ref{fig:some-paths-1} for an example, where the vertices from $U_1$ are drawn in black. 

	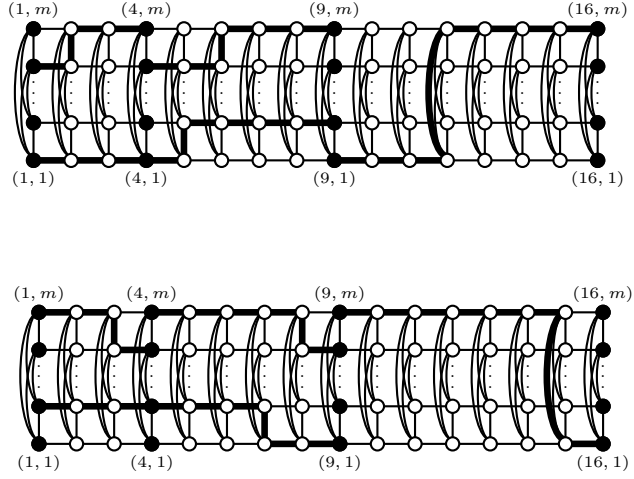
\begin{figure}[ht!]
		\centering
			\begin{tikzpicture}[thick,scale=0.5]
			\draw[line width = 2.5pt] (1,3.5) -- (2,3.5) -- (2,4.5) -- (4,4.5);
			\draw[line width = 2.5pt] (4,3.5) -- (6,3.5) -- (6,4.5) -- (9,4.5);
			\draw[line width = 2.5pt] (1,1) -- (5,1) -- (5,2) -- (9,2);
			\draw[line width = 2.5pt] (9,1) -- (12,1);
			\draw[line width = 2.5pt] (12,4.5) -- (16,4.5);
			\draw[line width = 2.5pt] (12,4.5) arc [x radius = 1.2em, y radius = 4.1em, start angle = 90, end angle = 270];
			\mreza{16}{0}
			\foreach \i in {1,2,3.5,4.5}
				{
				\draw [fill=black] (1,\i) circle(5pt);
				\draw [fill=black] (4,\i) circle(5pt);
				\draw [fill=black] (9,\i) circle(5pt);
				\draw [fill=black] (16,\i) circle(5pt);
				}
			\draw (1,0.5) node {\tiny $(1,1)$};
			\draw (4,0.5) node {\tiny $(4,1)$};
			\draw (9,0.5) node {\tiny $(9,1)$};
			\draw (16,0.5) node {\tiny $(16,1)$};
			\draw (1,5) node {\tiny $(1,m)$};
			\draw (4,5) node {\tiny $(4,m)$};
			\draw (9,5) node {\tiny $(9,m)$};
			\draw (16,5) node {\tiny $(16,m)$};
			\end{tikzpicture}
			\vspace{1cm}
			
			\begin{tikzpicture}[thick,scale=0.5]
			\draw[line width = 2.5pt] (1,4.5) -- (3,4.5) -- (3,3.5) -- (4,3.5);
			\draw[line width = 2.5pt] (4,4.5) -- (8,4.5) -- (8,3.5) -- (9,3.5);
			\draw[line width = 2.5pt] (1,2) -- (7,2) -- (7,1) -- (9,1);
			\draw[line width = 2.5pt] (9,4.5) -- (15,4.5);
			\draw[line width = 2.5pt] (15,1) -- (16,1);
			\draw[line width = 2.5pt] (15,4.5) arc [x radius = 1.2em, y radius = 4.1em, start angle = 90, end angle = 270];
			\mreza{16}{0}
			\foreach \i in {1,2,3.5,4.5}
				{
				\draw [fill=black] (1,\i) circle(5pt);
				\draw [fill=black] (4,\i) circle(5pt);
				\draw [fill=black] (9,\i) circle(5pt);
				\draw [fill=black] (16,\i) circle(5pt);
				}
			\draw (1,0.5) node {\tiny $(1,1)$};
			\draw (4,0.5) node {\tiny $(4,1)$};
			\draw (9,0.5) node {\tiny $(9,1)$};
			\draw (16,0.5) node {\tiny $(16,1)$};
			\draw (1,5) node {\tiny $(1,m)$};
			\draw (4,5) node {\tiny $(4,m)$};
			\draw (9,5) node {\tiny $(9,m)$};
			\draw (16,5) node {\tiny $(16,m)$};
			\end{tikzpicture}
			\\
		\caption{Shortest paths for $n=16$ and $(y_1,y_2,i,l)\in \{(m-1,m,2,1), (m-1,m,3,2), (1,2,3,1), (1,m,4,3)\}$.}
		\label{fig:some-paths-1}
	\end{figure}

\end{itemize}
The shortest paths from $\widehat{U_1}$ cover all the edges of $P_n \cp K_m$, hence we can conclude that $\sge(P_n \cp K_m) \le mk$ when $n = k^2$. 

Assume next that $n=k^2+h$, where $h \in [k]$. Then we claim that the set $U_2 = U_1 \cup \bigcup_{j=2}^m \{(n,j)\}$ is a strong edge geodetic set of $P_n \cp K_m$. 
\begin{itemize}
\item First put all the shortest paths from $\widehat{U_1}$ to $\widehat{U_2}$.
\item For every pair $y_1,y_2 \in V(K_m)$, where $y_1 < y_2$, and for every $i \in [h]$, put to $\widehat{U_2}$ the shortest path between the vertices $(i^2,y_1)$ and $(n,y_2)$ that contains the edge $(k^2+i,y_1)(k^2+i,y_2)$. See Fig.~\ref{fig:some-paths-2} for an example, where the vertices from $U_2$ are again drawn in black. 

	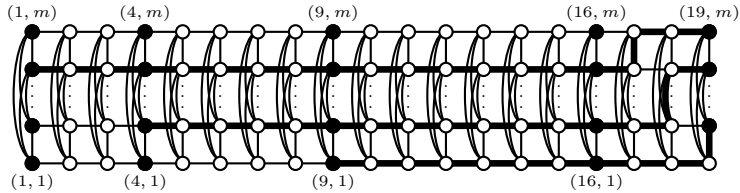
\begin{figure}[ht!]
		\centering
			\begin{tikzpicture}[thick,scale=0.5]
			\draw[line width = 2.5pt] (9,1) -- (19,1) -- (19,2);
			\draw[line width = 2.5pt] (1,3.5) -- (17,3.5) -- (17,4.5) -- (19,4.5);
			\draw[line width = 2.5pt] (4,2) -- (18,2);
			\draw[line width = 2.5pt] (18,3.5) -- (19,3.5);
			\draw[line width = 2.5pt] (18,3.5) arc [x radius = 0.4em, y radius = 1.7em, start angle = 90, end angle = 270]; 
			\mreza{19}{0}
			\foreach \i in {1,2,3.5,4.5}
				{
				\draw [fill=black] (1,\i) circle(5pt);
				\draw [fill=black] (4,\i) circle(5pt);
				\draw [fill=black] (9,\i) circle(5pt);
				\draw [fill=black] (16,\i) circle(5pt);
				}
			\draw [fill=black] (19,2) circle(5pt);
			\draw [fill=black] (19,3.5) circle(5pt);
			\draw [fill=black] (19,4.5) circle(5pt);
			\draw (1,0.5) node {\tiny $(1,1)$};
			\draw (4,0.5) node {\tiny $(4,1)$};
			\draw (9,0.5) node {\tiny $(9,1)$};
			\draw (16,0.5) node {\tiny $(16,1)$};
			\draw (1,5) node {\tiny $(1,m)$};
			\draw (4,5) node {\tiny $(4,m)$};
			\draw (9,5) node {\tiny $(9,m)$};
			\draw (16,5) node {\tiny $(16,m)$};
			\draw (19,5) node {\tiny $(19,m)$};
			\end{tikzpicture}
		\caption{Shortest paths for $n=19$ and $(y_1,y_2,i) \in \{(1,2,3), (2,m-1,2), (m-1,m,1)\}$.}
		\label{fig:some-paths-2}
	\end{figure}
	
\item For every $y \in \{2,\dots,m\}$, put to $\widehat{U_2}$ the unique shortest path between the vertices $(1,y)$ and $(n,y)$. Note that all the edges from $P_n^1$ are already covered by the shortest path from $\widehat{U_2}$ between vertices $(h^2,1)$ and $(n,2)$.
\end{itemize}

Since the shortest paths from $\widehat{U_2}$ cover all the edges of $P_n \cp K_m$, we can conclude that $\sge(P_n \cp K_m) \le mk + (m-1)$, when $n = k^2+h$ and $h \in [k]$.

Assume finally that $n=k^2+h$, where $k+1 \leq h \leq 2k$. In this case we claim that $U_3 = U_1 \cup \bigcup_{j=1}^m \{(n,j)\}$ is a strong edge geodetic set of $P_n \cp K_m$ and proceed as follows. 
\begin{itemize}
\item Put all the shortest paths from $\widehat{U_1}$ to $\widehat{U_3}$.
\item For every pair $y_1,y_2 \in V(K_m)$, where $y_1 < y_2$, do the following. For every $i \in [k]$, put to $\widehat{U_3}$ the shortest path between vertices $(i^2,y_1)$ and $(n,y_2)$ that contains the edge $(k^2+i,y_1)(k^2+i,y_2)$. Moreover, for every $i \in [h-k]$ also add the shortest path between the vertices $(n,y_1)$ and $(i^2,y_2)$ that contains the edge $(k(k+1)+i,y_1)(k(k+1)+i,y_2)$. In Fig.~\ref{fig:some-paths-3} examples are drawn with the vertices from $U_3$ again in black. 

	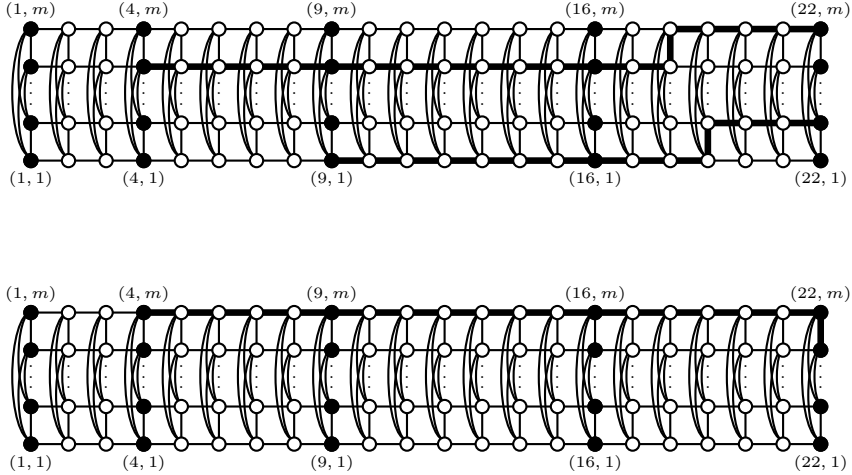
\begin{figure}[ht!]
		\centering
			\begin{tikzpicture}[thick,scale=0.5]
			\draw[line width = 2.5pt] (9,1) -- (19,1) -- (19,2) -- (22,2);
			\draw[line width = 2.5pt] (4,3.5) -- (18,3.5) -- (18,4.5) -- (22,4.5);
			\mreza{22}{0}
			\foreach \i in {1,2,3.5,4.5}
				{
				\draw [fill=black] (1,\i) circle(5pt);
				\draw [fill=black] (4,\i) circle(5pt);
				\draw [fill=black] (9,\i) circle(5pt);
				\draw [fill=black] (16,\i) circle(5pt);
				\draw [fill=black] (22,\i) circle(5pt);
				}
				\draw (1,0.5) node {\tiny $(1,1)$};
				\draw (4,0.5) node {\tiny $(4,1)$};
				\draw (9,0.5) node {\tiny $(9,1)$};
				\draw (16,0.5) node {\tiny $(16,1)$};
				\draw (22,0.5) node {\tiny $(22,1)$};
				\draw (1,5) node {\tiny $(1,m)$};
				\draw (4,5) node {\tiny $(4,m)$};
				\draw (9,5) node {\tiny $(9,m)$};
				\draw (16,5) node {\tiny $(16,m)$};
				\draw (22,5) node {\tiny $(22,m)$};
			\end{tikzpicture}
			
			\vspace{1cm}
			
			\begin{tikzpicture}[thick,scale=0.5]
			\draw[line width = 2.5pt] (4,4.5) -- (22,4.5) -- (22,3.5);
			\mreza{22}{0}
			\foreach \i in {1,2,3.5,4.5}
				{
				\draw [fill=black] (1,\i) circle(5pt);
				\draw [fill=black] (4,\i) circle(5pt);
				\draw [fill=black] (9,\i) circle(5pt);
				\draw [fill=black] (16,\i) circle(5pt);
				\draw [fill=black] (22,\i) circle(5pt);
				}
				\draw (1,0.5) node {\tiny $(1,1)$};
				\draw (4,0.5) node {\tiny $(4,1)$};
				\draw (9,0.5) node {\tiny $(9,1)$};
				\draw (16,0.5) node {\tiny $(16,1)$};
				\draw (22,0.5) node {\tiny $(22,1)$};
				\draw (1,5) node {\tiny $(1,m)$};
				\draw (4,5) node {\tiny $(4,m)$};
				\draw (9,5) node {\tiny $(9,m)$};
				\draw (16,5) node {\tiny $(16,m)$};
				\draw (22,5) node {\tiny $(22,m)$};
			\end{tikzpicture}
		\caption{Shortest paths for $n=22$ and $(y_1,y_2,i) \in \{(1,2,3), (m-1,m,2)\}$.}
		\label{fig:some-paths-3}
	\end{figure}
	
\item For every $j \in [m]$, put to $\widehat{U_3}$ the unique shortest path between the vertices $(1,j)$ and $(n,j)$. 
\end{itemize}
Since the shortest paths from $\widehat{U_3}$ cover all the edges of $P_n \cp K_m$, we get the upper bound $\sge(P_n \cp K_m) \le mk + m$ when $n = k^2+h$ with $k+1 \leq h \leq 2k$.

\medskip
In the second part of the proof we need to demonstrate that the obtained upper bounds are sharp, that is, there exist no smaller strong edge geodetic sets as the one constructed above. Let $U$ be a arbitrary strong edge geodetic set of $P_n \cp K_m$.

Assume first that $n = k^2$ for some $k \in \mathbb{N}$. Then we need to show that $|U| \geq m k$. If for every vertex $y \in K_m$, the set $U$ has at least $k$ vertices in the $P_n^y$-layer, then clearly $|U|\ge mk$. Assume therefore that for some $y_i \in V(K_m)$, the $P_n^{y_i}$- layer contains $k-l$, $l \geq 1$, vertices from $U$. Since $|V(P_n^{y_i}) \cap U| = k-l$, for every vertex $y \in V(K_m)$, $y\ne y_i$, the strong edge geodetic set $U$ has to have at least $x$ vertices from $P_n^{y}$, where $(k-l) x \geq k^2$, in order to cover all the edges between $P_n^{y_i}$ and $P_n^y$. Because $x \geq k^2/(k-l) = k + kl/(k-l) \geq k + kl/k = k+l$, we get 
$$|U| \ge k-l+(m-1)(k+l) = mk + (m-2)l \ge mk+1\,,$$ 
where the last assertion follows since $m \geq 3$ and $l \ge 1$. We conclude that in any case $|U| \ge mk$.  

Assume second that $n = k^2 + h$, where $1 \leq h \leq k$. Now we need to prove that $\sge(P_n \cp K_m) \geq m(k+1)-1$. If for every vertex $y \in K_m$, the set $U$ has at least $k+1$ vertices in $P_n^y$, then clearly $|U| \ge m(k+1)$ and we are done. Assume therefore that for some $y_i \in V(K_m)$, the set $U$ has $(k+1)-l$, $l \geq 1$, vertices from $P_n^{y_i}$. Since $|V(P_n^{y_i}) \cap U| \leq k+1-l$, for every vertex $y \in V(K_m)$, $y\ne y_i$, the set $U$ has to have at least $x$ vertices from $P_n^{y}$, where $(k+1-l)x \geq k^2+h$ has to hold in order to cover all the edges between $P_n^{y_i}$ and $P_n^y$. Because $x$ is an integer we can compute as follows:
\begin{align*}	
	x & \geq \left\lceil \frac{k^2+h}{k+1-l} \right\rceil
	 = \left\lceil \frac{k (k+h/k)}{k+1-l} \right\rceil 
	 = \left\lceil \frac{k (k+1-l) + k(h/k-1+l)}{k+1-l} \right\rceil \\
	 & = k + \left\lceil \frac{k(h/k-1+l)}{k+1-l} \right\rceil\,. 
\end{align*}
	Because $l \geq 1$ and therefore $1/(k+1-\ell) \geq 1/k$, we also have
$$x \geq k + \left\lceil \frac{k(h/k-1+l)}{k} \right\rceil = k + l-1+\left\lceil \frac{h}{k} \right\rceil\,.$$
	Since $h\in [k]$, we have $\left\lceil h/k \right\rceil = 1$ and therefore $x \geq k+l$. Altogether, 
\begin{align*}	
|U| & \ge k+1-l+(m-1)(k+l) = mk+(m-2)l+1 \\
& \geq mk+(m-2)+1 = m(k+1)-1
\end{align*}	
which we wanted to show. 
		
The remaining case is when $n = k^2 + h$, where $k+1 \leq h \leq 2k$. Now we need to prove that $\sge(P_n \cp K_m) \geq m(k+1)$. If for every vertex $y \in K_m$, the set $U$ has at least $k+1$ vertices from $P_n^y$,  then clearly $|U|\ge m(k+1)$. Assume therefore that for some $y_i \in V(K_m)$, the set $U$ has $(k+1)-l$, $l \geq 1$, vertices from $P_n^{y_i}$. So $|V(P_n^{y_i}) \cap U| \leq k+1-l$, hence for every vertex $y \in V(K_m)$, $y\ne y_i$, the set $U$ has to have at least $x$ vertices from  $P_n^{y}$, where $(k+1-l)x \geq k^2+h$ has to hold in order to cover all the edges between $P_n^{y_i}$ and $P_n^y$. Because $x$ is an integer and $l \geq 1$, we can similarly as in the previous case estimate that 
$$x \geq k + l-1+\left\lceil h/k \right\rceil\,.$$ 
Because $h$ is an integer between $k+1$ and $2k$, we have $\left\lceil h/k \right\rceil = 2$, and therefore $x \geq k + l + 1$. Altogether we see that 
\begin{align*}
|U| & \ge k+1-l+(m-1)(k+l+1)=m(k+1)+(m-2)l \\
& \geq m(k+1)+1\,,
\end{align*}
where the last assertion holds since $m \geq 3$ and $l \geq 1$. 
\qed

The following special case of Theorem~\ref{thm:path-by-complete} has been reported earlier in~\cite[Theorem~14]{xavier-2020}. 

\begin{corollary}
If $k\ge 2$ and $m\ge 3$, then $\sge(P_{k^2} \cp K_m) = mk$.
\end{corollary}

\section*{Acknowledgments} 

We would like to thank one of the reviewers for a very careful reading and many helpful tips, in particular on how to shorten the proof of Theorem~\ref{thm:complete multipartite}. 
This work was supported by the Slovenian Research Agency (ARRS) under the grants P1-0297, J1-2452, and N1-0285.

\section*{Declaration of interests}
 
The authors declare that they have no conflict of interest. 

\section*{Data availability}
 
Our manuscript has no associated data.

\end{document}